\documentclass[12pt,twoside, tikz, border=10pt]{article}
\usepackage{graphicx}
\usepackage{amssymb,array,arydshln}
\usepackage{amsfonts}
\usepackage{tikz}
\usetikzlibrary{chains,positioning}
\usepackage{amsbsy}
\usepackage{amssymb}\usepackage{color}
\usepackage{amsmath}
\usepackage{cite}
\usepackage{graphics}
\def\bpsp{\begin{pspicture}}
\def\epsp{\end{pspicture}}

\newtheorem{theorem}{Theorem}[section]
\newtheorem{remark}[theorem]{Remark}
\newtheorem{example}[theorem]{Example}
\newtheorem{lemma}[theorem]{Lemma}
\newtheorem{corollary}[theorem]{Corollary}
\newtheorem{definition}[theorem]{Definition}
\newtheorem{proposition}[theorem]{Proposition}
\newtheorem{note}{Note}
\newtheorem{case}{Case}
\newtheorem{conjecture}{Conjecture}
\newtheorem{question}{Question}

\newcommand{\bea}{\begin{eqnarray}}
\newcommand{\eea}{\end{eqnarray}}
\newcommand{\beq}{\begin{eqnarray*}}
\newcommand{\eeq}{\end{eqnarray*}}

\def\m4{\mbox{\rm ~(mod $4$)}}

\def \bd{\begin{definition}}
\def \ed{\end{definition}}
\def \bqu{\begin{question}}
\def \equ{\end{question}}
\def \bcc{\begin{conjecture}}
\def \ecc{\end{conjecture}}
\def \bt{\begin{theorem}}
\def \et{\end{theorem}}
\def \bl{\begin{lemma}}
\def \el{\end{lemma}}
\def \bc{\begin{corollary}}
\def \ec{\end{corollary}}
\def \be{\begin{equation}}
\def \ee{\end{equation}}
\def \ben{\begin{enumerate}}
\def \een{\end{enumerate}}
\def \ba{\begin{array}}
\def \ea{\end{array}}
\def \bp{\begin{proposition}}
\def \ep{\end{proposition}}
\def \bx{\begin{example}}
\def \ex{\end{example}}
\def \br{\begin{remark}}
\def \er{\end{remark}}
\def \bdsc{\begin{description}}
\def \edsc{\end{description}}

\def \bn{\begin{case}}
\def \en{\end{case}}
\def \bnt{\begin{note}}
\def \ent{\end{note}}
\def\1{1\!\!1}

\def\mm2{\mbox{\rm ~(mod $2$)}}
\def\m4{\mbox{\rm ~(mod $4$)}}

\def\qed{\nolinebreak\hfill\rule{.2cm}{.2cm}\par\addvspace{.5cm}}

\def\m{\mu}

\def\1{\textbf{1}}
\def\0{\textbf{0}}

\begin{document}

\title{Direct sum graph of the subspaces of a finite dimensional vector space over finite fields}
\author{Bilal A. Wani$^a$, Aaqib Altaf$^b$, S. Pirzada$^c$, T. A. Chishti$^d$\\
$^{a}${\em Department of Mathematics, National Institute of Technology Srinagar,} \\
{\em Srinagar, Kashmir, India}\\
	$^{b,c,d}${\em Department of Mathematics, University of Kashmir, Srinagar, India}\\
	$^a$bilalwanikmr@gmail.com; $^b$aaqibwani777@gmail.com\\
$^c$pirzadasd@kashmiruniversity.ac.in; $^d$tachishti@uok.edu.in}
\date{}

\pagestyle{myheadings} \markboth{Wani, Altaf, Pirzada, Chishti}{Direct sum graph of the subspaces of a finite dimensional vector space}
\maketitle
\vskip 5mm
\noindent{\footnotesize \bf Abstract.} In this paper, we introduce a new graph structure, called the $direct~ sum ~graph$ on a finite dimensional vector space. We investigate the connectivity, diameter
and the completeness of $\Gamma_{U\oplus W}(\mathbb{V})$. Further, we find its domination number and independence
number. We also determine the degree of each vertex in case the base field is finite and show that the graph $\Gamma_{U\oplus W}(\mathbb{V})$ is  not Eulerian. We also show that under some mild conditions the graph  $\Gamma_{U\oplus W}(\mathbb{V})$ is triangulated. We determine the clique number of $\Gamma_{U\oplus W}(\mathbb{V})$ for some particular cases. Finally, we find the size, girth, edge-connectivity and the chromatic number of $\Gamma_{U\oplus W}(\mathbb{V})$.
\vskip 3mm

\noindent{\footnotesize Keywords: Vector space; graph, Maximal cliques;  Euler graph.}

\vskip 3mm
\noindent {\footnotesize AMS subject classification: 05C25, 05C69.}

\section{Introduction}
Let $G=(V(G), E(G))$ be a simple graph with vertex set $V(G)=\{v_{1},v_{2},\ldots,v_{n}\}$ and edge set $E(G)$. For $1\leq i\neq j\leq n$, we write $v_i\thicksim v_j$ if $v_i$ is adjacent to $v_j$ in $G$. The cardinality of $V(G)$ and  $E(G)$  is called the \textit{order} and the \textit{size} of $G$,  respectively. If $u\in V(G)$, then $N(u)$ is the set of neighbors of $u$ in $G$, that is, $N(u)=\{v\in (G): uv\in E(G)\}$.  For any two distinct vertices $u$ and $v$ of $G$,  $d(u,v)$ denotes the length of a shortest path between $u$ and $v$. Clearly $d(u,u)=0$ and $d(u,v)=\infty$, if there is no path connecting $u$ and $v$. The $diameter$ of $G$ is defined as
$diam(G)=max \{d(u,v): u,v \in V(G)\}$.  A graph $G$ is said to be $complete$ if every pair of distinct vertices are adjacent and a complete graph on $n$ vertices is denoted by $K_n$. A complete subgraph of a graph $G$ is called a clique. A maximal clique is a clique
which is maximal with respect to inclusion. The clique number of $G$, written as $\omega(G)$,
is the maximum order of a clique in $G$. The chromatic number of $G$, denoted as $\chi(G)$,
is the minimum number of colors required to label the vertices so that the adjacent
vertices receive different colors. A graph $G$ is said to be $connected$ if for any pair of vertices $u, v \in V(G)$, there exists a path between $u$ and $v$. A subset $\alpha(G)$ of $V$ is said to be independent if no two vertices in that subset are pairwise adjacent. The independence number of a graph is the cardinality of the maximum independent set of the vertices in $G$. A subset $D$ of $V$ is said to be dominating set if every edge in $V\setminus D$ is adjacent to at least one vertex in $D$. The dominating number of $G$, denoted by $\gamma(G)$, is the cardinality of the minimum dominating set in $G$. A subset $D$ of $V$ is said to be a minimal dominating set if $D$ is a dominating set and no proper subset of $D$ is a dominating set.  A graph is said to be Eulerian if it contains a closed walk which traverses all the edges in $G$ exactly once. A graph is said to be triangulated if for any vertex $u$ in $V$, there exists $v, w\in V$, such that $(u,v,w)$ is a triangle. The girth of a graph is the length of the shortest cycle, if it exists. Otherwise, it is defined as $\infty$. For undefined terms and concepts, the reader is referred to \cite{GP,sp}.

\indent In 1988, Beck \cite{beck} initiated the study of graphs associated with  various algebraic structures by introducing the concept of zero-divisor graphs associated to rings. Redmond \cite{red2} extended the study of the zero divisor graph of commutative rings to non-commutative rings. DeMeyer et al. \cite{DMS} associated graphs to semigroups. Redmond \cite{red3} took a new approach to define zero divisor graphs with the help of an ideal of a ring. The association of a graph and vector space has history back in 1958 by Gould \cite{GOU}. Later, Chen \cite{CHE} investigated on vector spaces associated with a graph. Carvalho \cite{CAR} has studied vector space and the Petersen Graph. In the
recent past, Manjula \cite{MAN} used vector spaces and made it possible to use techniques
of linear algebra in studying the graph. Intersection graphs associated with subspaces of vector spaces were studied in \cite{CAGH,CHGHMS1,CHGHMS2, JRJS,TEA}. Recently, Das \cite{DAS01} introduced the graphs associated with the elements of finite dimensional vector space known as $nonzero~ component~graph$ of vector space. More work on this can be seen in \cite{DAS03, DAS04,DAS05}.

Motivated by the above work, we introduce the $direct~sum~graph$ as follows: Let $\mathbb{V}$ be a vector space over a field $\mathcal{F}$ with two subspaces $U$ and $W$. We say that $\mathbb{V}$ is the direct sum of the subspaces $U$ and $W$, if each vector of $\mathbb{V}$ can always be expressed as the sum of a vector from $U$ and a vector from $W$, and this expression can only be accomplished in one way (i.e. uniquely). In other words, for every $x\in \mathbb{V}$, there exists vectors $u\in U$, and $w\in W$ such that $x=u+w$ and this representation of $v$ is unique. Let $\mathcal{B_\mathbb{V}}=\left\lbrace \alpha_1, \alpha_2, \dots, \alpha_r+\beta_1+\beta_2+\dots+\beta_s\right\rbrace $, with $r+s=n$ be a basis of $\mathbb{V}$ and $\mathcal{B_\mathbb{U}}=\left\lbrace \alpha_1, \alpha_2, \dots, \alpha_r\right\rbrace $ and $\mathcal{B_\mathbb{W}}=\left\lbrace \beta_{1}, \beta_{2}, \dots, \beta_s\right\rbrace$ be basis of $U$ and $W$ respectively. Then any vector $x=u+w$ with unique representation as $\mathbf{x}=a_1\alpha_1+a_2\alpha_2+\dots+a_r\alpha_r+b_{1}\beta_{1}+b_{2}\beta_{2}+\dots+b_s\beta_s$, is said to have its basic representation with respect to $\mathcal{B_\mathbb{U}}$ and $\mathcal{B_\mathbb{W}}$.  We define the $direct~ sum ~graph$, denoted by $\Gamma_{U\oplus W}(\mathbb{V})=\left( V,E\right)$, of a finite dimensional vector space  with respect to $U$ and  $W$ as follows: $V=\{\mathbf{x}=u+w\in \mathbb{V}~| u\neq 0~ \text{and}~w\neq 0\}$ and for $\mathbf{x_1},\mathbf{x_2}\in V$, there is an edge between $\mathbf{x_1}$ and $\mathbf{x_2}$, that is, $\mathbf{x_1}\sim \mathbf{x_2}$ or $\left( \mathbf{x_1},\mathbf{x_2}\right) \in E$ if and only if $\mathbf{x_1}$ and $\mathbf{x_2}$ share at least two basis elements one each from  $\mathcal{B_\mathbb{U}}$  and $\mathcal{B_\mathbb{W}}$ having  nonzero  coefficients in there basic representation. It is necessary to mention here that one of the subspaces is the zero subspace, that is, $\mathbb{U}={0}$(say),  that is, $\mathbf{x_1}\sim \mathbf{x_2}$ if and only if $\mathbf{x_1}$ and $\mathbf{x_2}$ share at least one basis element  from $\mathcal{B_\mathbb{W}}$ having  nonzero  coefficients in there basic representation. In other words, the  $direct~ sum ~graph$ is the generalization of the $nonzero~  component ~ graph$.   For some examples of the $direct~ sum ~graph$  $\Gamma_{U\oplus W}(\mathbb{V})$, see Figures 1 and 2.

\begin{figure}
	\centering
	\begin{tikzpicture}
	\draw[fill=black] (0,0) circle (2pt);\node at (0.3,0.3) {$\alpha_{1}+\beta_{1}$};
	
	\draw[fill=black] (3,0) circle (2pt);\node at (3.2,0.3) {$\alpha_{2}+\beta_{1}$};
	
	\draw[fill=black] (7,0) circle (2pt); \node at (7.3, 0.3) {$\alpha_{1}+\beta_{2}$};
	\draw[fill=black] (11,0) circle (2pt); \node at (11.3,0.3) {$\alpha_{2}+\beta_{2}$};

		\draw[fill=black] (1.5, 2) circle (2pt);\node at (1.8,2.3) {$\alpha_{1}+\alpha_{2}+\beta_{1}$};
	
	\draw[fill=black] (8.5,2) circle (2pt);\node at (8.8,2.3) {$\alpha_{1}+\alpha_{2}+\beta_{2}$};
	\draw[fill=black] (1.5,4) circle (2pt);\node at (1.8,4.3) {$\alpha_{1}+\beta_{1}+\beta_{2}$};

	\draw[fill=black] (8.5,4) circle (2pt); \node at (8.8, 4.3) {$\alpha_{2}+\beta_{1}+\beta_{2}$};
	\draw[fill=black] (4.5,7) circle (2pt); \node at (4.8,7.3) {$\alpha_{1}+\alpha_{2}+\beta_{1}+\beta_{2}$};
	
	\draw[thin] (4.5,7)--(8.5,4); \draw[thin] (4.5,7)--(1.5,4);
	\draw[thin] (4.5,7)--(8.5,2); \draw[thin] (4.5,7)--(1.5,2);
	\draw[thin] (4.5,7)--(11,0); \draw[thin] (4.5,7)--(7,0);
	\draw[thin] (4.5,7)--(3,0); \draw[thin] (4.5,7)--(0,0);
	\draw[thin] (1.5,4)--(1.5,2); \draw[thin] (1.5,4)--(7,0);
	\draw[thin] (1.5,4)--(0,0); \draw[thin] (8.5,4)--(8.5,2);
	\draw[thin] (8.5,4)--(3,0); \draw[thin] (8.5,4)--(11,0);
	\draw[thin] (1.5,2)--(0,0); \draw[thin] (1.5,2)--(3,0);
	\draw[thin] (8.5,2)--(7,0); \draw[thin] (8.5,2)--(11,0);
	\end{tikzpicture}
	\caption{$dim(\mathbb{V})=4,~ ,~dim(U)=2,~ dim(W)=2,~ and ~ \mathcal{F}= \mathbb{F}_2$}
\end{figure}
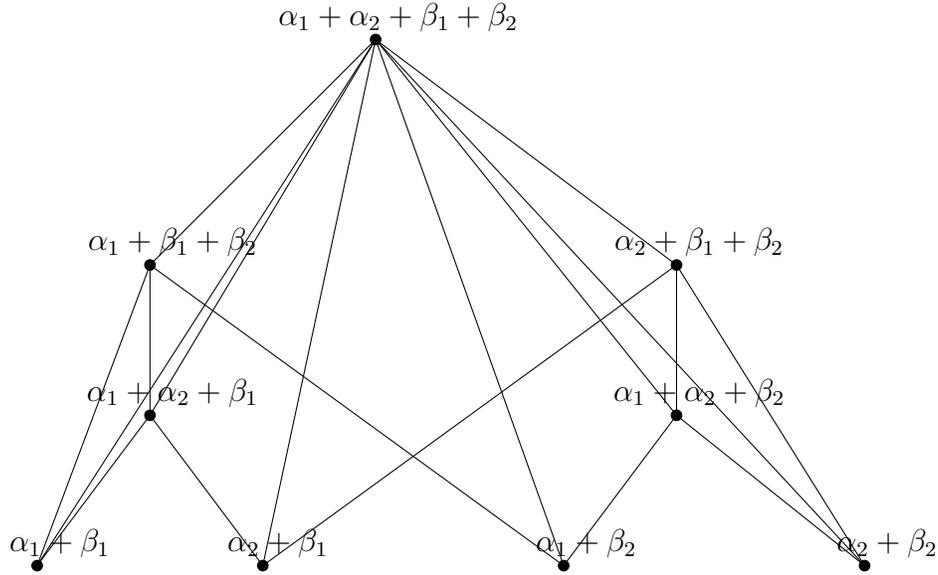
\begin{figure}
	\centering
	\begin{tikzpicture}
	\draw[fill=black] (1.5,2) circle (2pt); \node at (1.5, 2.3) {$\alpha_{1}+\beta_{1}$};
	\draw[fill=black] (4.5,2) circle (2pt); \node at (4.5,2.3) {$\alpha_{1}+\beta_{2}$};
	\draw[fill=black] (7.5,2) circle (2pt); \node at (7.5,2.3) {$\alpha_{1}+\beta_{3}$};
	\draw[fill=black] (0,0) circle (2pt);\node at (-1.4,0) {$\alpha_{1}+\beta_{1}+\beta_{2}$};
	\draw[fill=black] (3.5,0) circle (2pt);\node at (3.3,0.3) {$\alpha_{1}+\beta_{1}+\beta_{3}$};
	\draw[fill=black] (6,0) circle (2pt);\node at (7.2,0) {$\alpha_{1}+\beta_{2}+\beta_{3}$};
	\draw[fill=black] (1.6,-2) circle (2pt); \node at (1.5,-2.3) {$\alpha_{1}+\beta_{1}+\beta_{2}+\beta_{3}$};
	\draw[thin] (1.5,2)--(0,0); \draw[thin] (1.5,2)--(3.5,0);
	\draw[thin] (1.6,-2)--(1.5,2); \draw[thin] (4.5,2)--(0,0);
	\draw[thin] (1.6,-2)--(4.5,2); \draw[thin] (4.5,2)--(6,0);
	\draw[thin] (7.5,2)--(3.5,0); \draw[thin] (7.5,2)--(6,0);
	\draw[thin] (7.5,2)--(1.6,-2); \draw[thin] (1.6,-2)--(0,0);
	\draw[thin] (1.6,-2)--(3.5,0);\draw[thin] (1.6,-2)--(6,0);
	\draw[thin] (0,0)--(3.5,0);\draw[thin] (0,0)--(6,0);
	\draw[thin] (3.5,0)--(6,0);
	\end{tikzpicture}
	\caption{$dim(\mathbb{V})=3,~dim(U)=1,~ dim(W)=2,~ and ~ \mathcal{F}= \mathbb{F}_2$}
\end{figure}
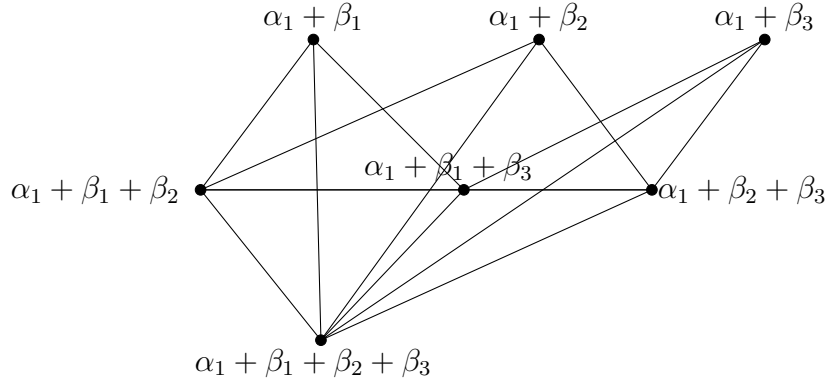
\indent Throughout this paper, the subspaces $U$ and $W$ of $\mathbb{V}=U\oplus W$ are non trivial and $dim(\mathbb{V})>1$.

The paper is organized as follows. In Section 2, we determine the  order, size and  edge connectivity of $\Gamma_{U\oplus W}(\mathbb{V})$.  In Section 3, we show that $\Gamma_{U\oplus W}(\mathbb{V})$ is  not Eulerian. In addition, we find the conditions under which the graph $\Gamma_{U\oplus W}(\mathbb{V})$ is triangulated. In Section 4, we find clique  number and the chromatic number of $\Gamma_{U\oplus W}(\mathbb{V})$.

\section{Properties of $\Gamma_{U\oplus W}(\mathbb{V})$}

We begin this section by investigating some basic properties like connectedness, com-pleteness, independence number and domination number of $\Gamma_{U\oplus W}(\mathbb{V})$.
\begin{theorem}\label{T1} Let $\mathbb{V}=U\oplus W$ be a vector space over a field $\mathcal{F}$. Then $\Gamma_{U\oplus W}(\mathbb{V})$ is connected and $diam(\Gamma_{U\oplus W})=2$, if $dim(\mathbb{V})\geq3$. If $dim(\mathbb{V})=2$, then $diam(\Gamma_{U\oplus W})=1.$
\end{theorem}
\noindent\textbf{Proof.} Assume that $dim(\mathbb{V})=2,$ $U$ and $W$ are non trivial subspaces of $\mathbb{V}$ that is $dim(U)=dim(W)=1$. It is easy to see that any two vertices $\mathbf{x}, \mathbf{y}\in \mathbf{V}$ are adjacent. Hence, $\Gamma_{U\oplus W}(\mathbb{V})$ is a complete graph. Now, assume that $dim(\mathbb{V})=n\geq 3$.   Let $Dim(U)=r$ and  $Dim(W)=s$, with $r,s\geq 1$ and $r+s=n$.  Depending on the choice of a subspace $U$ and $W$, we will break the proof in two cases:\\
\textbf{Case(i)} Assume that $Dim(U)=1$. Then $Dim(W)=n-1$, that is, $W$ is a hyperspace. Let  $\mathbf{x}, \mathbf{y}\in \mathbf{V}$ such that $\mathbf{x}=u+v$ and $\mathbf{y}=u'+v'$. If $\mathbf{x}$ and $\mathbf{y}$ are adjacent,  then $d(\mathbf{x}, \mathbf{y})=1.$ Let $\mathbf{x}$ and $\mathbf{y}$ bee not adjacent. Since $dim(U)=1$, so $u$ and $u'$ share non-zero coefficient in their basic representation, but $v$ and $v'$ do not share any  non-zero coefficient in their basic representation, that is, there exists $\beta_i, \beta_j$ ($\beta_i\neq \beta_j$), which have non-zero coefficient in the basic representation of $v$ and $v'$, respectively. Consider $\mathbf{z}=\mathbf{x}+\mathbf{y}=u+u'+v+v'$. Clearly,  $\mathbf{x}\sim \mathbf{z}$ and $\mathbf{z}\sim \mathbf{y}$. Hence, $d(\mathbf{x}, \mathbf{y})=2.$ Thus, $\Gamma_{U\oplus W}(\mathbb{V})$ is connected and $diam(\Gamma_W)=2$.\\
\textbf{Case(ii)} Assume that $Dim(U)=r$ and  $Dim(W)=s$, with $r,s\geq 2$ and $r+s=n$. Let  $\mathbf{x}, \mathbf{y}\in \mathbf{V}$ such that $\mathbf{x}=u+v$ and $\mathbf{y}=u'+v'$. If $\mathbf{x}$ and $\mathbf{y}$ are adjacent,  then $d(\mathbf{x}, \mathbf{y})=1.$ If $\mathbf{x}$ and $\mathbf{y}$ are not adjacent, then both $u,u'$ and $v,v'$ do not share a non-zero coefficient in their basic representation, that is, there exists $\alpha_i,\alpha_j$ and $\beta_i, \beta_j$ ($\alpha_i\neq \alpha_j$ and  $\beta_i\neq \beta_j$), which have non-zero coefficient in the basic representation of $u,u'$ and $v, v'$, respectively. Consider $\mathbf{z}=\mathbf{x}+\mathbf{y}=u+u'+v+v'$. Clearly,  $\mathbf{x}\sim \mathbf{z}$ and $\mathbf{z}\sim \mathbf{y}$. Hence, $d(\mathbf{x}, \mathbf{y})=2.$ Thus, $\Gamma_{U\oplus W}(\mathbb{V})$ is connected and $diam(\Gamma_W)=2$.\qed

\begin{theorem}\label{T2} Let $\mathbb{V}=U\oplus W$ be a vector space over a field $\mathcal{F}$. Then $\Gamma_{U\oplus W}(\mathbb{V})$  is complete if and only if $dim(\mathbb{V})=2$.
\end{theorem}
\noindent\textbf{Proof.} Assume that $\Gamma_{U\oplus W}(\mathbb{V})$ is complete. We will show that  $dim(\mathbb{V})=2$. Suppose on the contrary that $dim(\mathbb{V})>2$. Let $Dim(U)=r$ and  $Dim(W)=s$, with $r\geq 1$, $s\geq 2$ and $r+s=n$. Let  $\mathbf{x}, \mathbf{y}\in \mathbf{V}$ such that $\mathbf{x}=\alpha_i+\beta_j$ and $\mathbf{y}=\alpha_i+\beta_k$, where $j\neq k$. Clearly, $\mathbf{x}\nsim \mathbf{z}$, a contradiction. Thus,   $dim(\mathbb{V})=2$.

Conversely, suppose that $dim(\mathbb{V})=2$. So $U$ and $W$ are non trivial subspaces of $\mathbb{V}$. In other words, $dim(U)=dim(W)=1$. Any two vertices $\mathbf{x}$ and $\mathbf{y}$ are of the form $\mathbf{x}=a_i\alpha_1+b_i\beta_1$ and  $\mathbf{y}=a'_i\alpha_1+b'_i\beta_1$. Therefore $\mathbf{x}\sim \mathbf{y}$. Hence, $\Gamma_{U\oplus W}(\mathbb{V})$ is a complete graph.\qed
\begin{theorem}\label{T3} Let $\mathbb{V}=U\oplus W$ be a vector space over a field $\mathcal{F}$. Then the domination number of $\Gamma_{U\oplus W}(\mathbb{V})$ is 1.
\end{theorem}
\noindent\textbf{Proof.} Assume that $Dim(U)=r$ and  $Dim(W)=s$, with $r+s=n$. Let $\{\alpha_1,\alpha_2,\dots,\alpha_r\}$ and $\{\beta_1,\beta_2,\dots,\beta_s\}$ be the basis of $U$ and $W$, respectively. It is  easy to see that the vertex $\alpha_1+\alpha_2+\dots+\alpha_r+\beta_1+\beta_2+\dots+\beta_s$ is adjacent to any other vertex of $\mathbf{V}$. Thus the domination number of $\Gamma_{U\oplus W}(\mathbb{V})$ is 1.\qed
\begin{remark} The set $\{\alpha_1+\beta_1,\alpha_1+\beta_2,\dots, \alpha_1+\beta_r,\alpha_2+\beta_1,\alpha_2+\beta_2,\dots,\alpha_r+\beta_s\}$ is a minimal dominating set of $\Gamma_{U\oplus W}(\mathbb{V})$. Now, the
question arises what is the maximum possible number of vertices in a minimal dominating
set. The answer is given as $dim(U)dim(W)$ in the next theorem.
\end{remark}
\begin{theorem}\label{T4} Let $\mathbb{V}=U\oplus W$ be a vector space over a field $\mathcal{F}$. If $\mathbf{D}=\{x_1,x_2,\dots,x_k\}$ is a minimal dominating set of  $\Gamma_{U\oplus W}(\mathbb{V})$, then $k \leq dim(U)dim(W)$.
\end{theorem}
\noindent\textbf{Proof.} Assume that $dim(U)=r$ and  $dim(W)=s$, with $r+s=n$. Let $\{\alpha_1,\alpha_2,\dots,\alpha_r\}$ and $\{\beta_1,\beta_2,\dots,\beta_s\}$ be the basis of $U$ and $W$, respectively. Let $\mathbf{D}=\{x_1,x_2,\dots,x_k\}$ be a minimal dominating set of  $\Gamma_{U\oplus W}(\mathbb{V})$. Construct $\mathbf{D}_i=\mathbf{D}\setminus\{x_i\}$, for $1\leq i\leq k$. Indeed,   $\mathbf{D}_i$ is not a dominating set. In other words, corresponding to every $\mathbf{D}_i$, there exists $\eta_i\in \Gamma_{U\oplus W}(\mathbb{V})$, which is not adjacent to any
element of $\mathbf{D}_i$ but adjacent to $x_i$. Since, $\eta_i\neq 0$, there exit $\alpha^i_{p}$ and $\beta^i_{q}$ such that $\eta_i$ has non-zero component along $\alpha^i_{p}$ and $\beta^i_{q}$. Now, as $\eta_i$ is not adjacent to any element of $\mathbf{D}_i$, so is $\alpha^i_{p}+\beta^i_{q}$.  Thus,  for all $i\in \{1, 2, \dots, k\}$, there exits $\alpha^i_{p}+\beta^i_{q}$
such that $\alpha^i_{p}+\beta^i_{q}\sim x_i$, but $\alpha^i_{p}+\beta^i_{q}\nsim x_j$, for all $j\neq i$.

 To complete the proof, it only suffices to show that $\alpha^i_{p}+\beta^i_{q}\neq \alpha^j_{p}+\beta^j_{q}$ for $i\neq j$. Suppose on contrary $\alpha^i_{p}+\beta^i_{q}= \alpha^j_{p}+\beta^j_{q}$ for $i\neq j$. Since
$\alpha^i_{p}+\beta^i_{q}\sim x_i$ and $\alpha^i_{p}+\beta^i_{q}= \alpha^j_{p}+\beta^j_{q}$, it follows that $\alpha^j_{p}+\beta^j_{q}\sim x_i$, which contradicts $\alpha^i_{p}+\beta^i_{q}\nsim x_j$, for all $j\neq i$. Hence, $\alpha^i_{p}+\beta^i_{q}\neq \alpha^j_{p}+\beta^j_{q}$ for $i\neq j$. i.e, all $\alpha^i_{p}+\beta^i_{q}$'s are distinct, and the total number of such combinations are $rs$, it follows that $k\leq dim(U)dim(W)$.\qed
From Theorem \ref{T2}, it is clear that $\Gamma_{U\oplus W}(\mathbb{V})$  is complete if and only if $dim(\mathbb{V})=2$. Now, the question arises what is the independence number  of $\Gamma_{U\oplus W}(\mathbb{V})$ if $dim(\mathbb{V})\geq 3$. The answer is given in the next theorem.

\begin{theorem}\label{T5} Let $\mathbb{V}=U\oplus W$ be a vector space over a field $\mathcal{F}$ with $dim(\mathbb{V})\geq3$. Then the independence number of  $\Gamma_{U\oplus W}(\mathbb{V})$ is $dim(U)dim(W)$.
\end{theorem}
\noindent\textbf{Proof.} Since $\{\alpha_1+\beta_1,\alpha_1+\beta_2,\dots, \alpha_1+\beta_r,\alpha_2+\beta_1,\alpha_2+\beta_2,\dots,\alpha_r+\beta_s\}$ is an independent set in $\Gamma_{U\oplus W}(\mathbb{V})$, the independence number of $\Gamma_{U\oplus W}(\mathbb{V})\geq dim(U)dim(W)$.  It suffices to show that the independence number of $\Gamma_{U\oplus W}(\mathbb{V})\leq dim(U)dim(W)$. Suppose to the contrary that $\{ x_1, x_2, \dots, x_k\}$ is an independent set in $\Gamma_{U\oplus W}(\mathbb{V})$ with $k>dim(U)dim(W)$, that is, $x_i \nsim x_j$ for $1\leq i\neq j \leq k$. Also, for $x_i\in \mathbf{V}$,  for all $i\in \{1,2,\dots,k\}$, there exits $\alpha^i_{p}, \beta^i_{q}$ having non zero coefficient along $\alpha^i_{p}$ and $\beta^i_{q}$ in its basic representation.  Claim that $\alpha^i_{p}+\beta^i_{q}\neq \alpha^j_{p}+\beta^j_{q}$ for $i\neq j$. Assume that $\alpha^i_{p}+\beta^i_{q}= \alpha^j_{p}+\beta^j_{q}$ for $i\neq j$. Since
$\alpha^i_{p}+\beta^i_{q}\sim x_i$ and $\alpha^i_{p}+\beta^i_{q}= \alpha^j_{p}+\beta^j_{q}$, it follows that $\alpha^j_{p}+\beta^j_{q}\sim x_i$. But $\alpha^j_{p}+\beta^j_{q}\sim x_j$, which contradicts $x_i\nsim x_j$, for all $i\neq j$. Hence, $\alpha^i_{p}+\beta^i_{q}\neq \alpha^j_{p}+\beta^j_{q}$ for $i\neq j$.   As there are exactly  $dim(U)dim(W)$ distinct $\alpha_i+\beta_j$'s, it follows that $k\leq dim(U)dim(W)$, which is a contradiction. Hence, the independence number of  $\Gamma_{U\oplus W}(\mathbb{V})$ is $dim(U)dim(W)$.\qed
Now, we will find the degree of each vertex of $\Gamma_{U\oplus W}(\mathbb{V})$ if the base field is finite.  Let  $\mathcal{F}=\mathbb{F}_q$ be a field  with $q$ elements and $\mathbb{V}=U\oplus W$  be an $n$ dimensional vector space over $\mathcal{F}$. Let $dim(U)=r$ and  $dim(W)=s$, with $r+s=n$. Let $\{\alpha_1,\alpha_2,\dots,\alpha_r\}$ and $\{\beta_1,\beta_2,\dots,\beta_s\}$ be basis of $U$ and $W$, respectively.
It is to be noted that
\begin{eqnarray}
\nonumber \lefteqn{ N(\alpha_{i_1}+\alpha_{i_2}+\dots+ \alpha_{i_l}+\beta_{j_1}+\beta_{j_2}+\dots+ \beta_{j_m})}\\
\nonumber &=& N(c_1\alpha_{i_1}+c_{2}\alpha_{i_2}+\dots+ c_{l}\alpha_{i_l}+d_{1}\beta_{j_1}+d_{2}\beta_{j_2}+\dots+ d_{m}\beta_{j_m}),
\end{eqnarray}
where $1\leq l\leq r$ and $1\leq m\leq s$ and  $c_i\neq0\neq d_j\in \mathbb{F}.$

\begin{theorem}\label{T7} Let $\mathbb{V}=U\oplus W$ be an $n$ dimensional vector space over a field $\mathcal{F}$ with $q$ elements. Suppose  $\{\alpha_1,\alpha_2,\dots,\alpha_r\}$ and $\{\beta_1,\beta_2,\dots,\beta_s\}$ are the basis of $U$ and $W$, respectively, with $r+s=n$. Let $\Gamma_{U\oplus W}(\mathbb{V})$  be
a graph  associated with $\mathbb{V}=U\oplus W$ Then, the degree of the vertex $c_1\alpha_{i_1}+c_{2}\alpha_{i_2}+\dots+ c_{l}\alpha_{i_l}+d_{1}\beta_{j_1}+d_{2}\beta_{j_2}+\dots+ d_{m}\beta_{j_m}$, where $c_1c_2\dots c_l d_1d_2\dots d_m \neq 0$, is $(q^l-1)(q^m-1)q^{n-(l+m)}-1$.
\end{theorem}
\noindent\textbf{Proof.} First fix $\beta_{j_1}$. The number of vertices with $\alpha_{i_1}$ and $\beta_{j_1}$, as a nonzero component in the basic representation (including $\alpha_{i_1}+\beta_{j_1}$
 itself) is $(q-1)^2q^{r+s-2}$. Therefore, $$deg(\alpha_{i_1}+\beta_{j_1})=(q-1)^2q^{r+s-2}-1.$$
The number of the vertices with ($\alpha_{i_1}$ and $\beta_{j_1}$) or ($\alpha_{i_2}$ and $\beta_{j_1}$) as the non-zero component is equal to the number
of vertices with ($\alpha_{i_1}$ and $\beta_{j_1}$) as the non-zero component $+$ the number of vertices with ($\alpha_{i_2}$ and $\beta_{j_1}$) as the non-zero
component $-$ the number of vertices with both ($\alpha_{i_1}, \alpha_{i_2}$ and $\beta_{j_1}$)  as the non-zero component
\begin{align*}
&= (q-1)^2q^{r+s-2}+(q-1)^2q^{r+s-2}-(q-1)^3q^{r+s-3}\\&
=2(q-1)^2q^{r+s-2}-(q-1)^3q^{r+s-3}\\&
=\big(2(q-1)q^{r-1}-(q-1)^2q^{r-2}\big)(q-1)q^{s-1}\\&
= (q^2-1)q^{r-2}(q-1)q^{s-1}.
\end{align*}
Therefore,
$deg(\alpha_{i_1}+ \alpha_{i_2}+\beta_{j_1})=(q^2-1)q^{r-2}(q-1)q^{s-1}-1.$ \\
Similarly, for finding the degree of $\alpha_{i_1}+ \alpha_{i_2}+\alpha_{i_3}+\beta_{j_1}$, the number of vertices with $\alpha_{i_1}\text{and} \beta_{j_1}$ or
$\alpha_{i_2}\text{and} \beta_{j_1}$ or $\alpha_{i_3}\text{and} \beta_{j_1}$ as the non-zero component is equal to
\begin{align*}
&=\binom{3}{1}(q-1)^2q^{r+s-2}-\binom{3}{2}(q-1)^3q^{r+s-3}+\binom{3}{3}(q-1)^4q^{r+s-4}\\&
=\Bigg[\binom{3}{1}(q-1)q^{r-1}-\binom{3}{2}(q-1)^2q^{r-2}+\binom{3}{3}(q-1)^3q^{r-3}\bigg](q-1)q^{s-1}\\&
=(q^3-1)q^{r-3}(q-1)q^{s-1}.
\end{align*}
Hence, $deg(\alpha_{i_1}+ \alpha_{i_2}+\alpha_{i_3}+\beta_{j_1})=(q^3-1)q^{r-3}(q-1)q^{s-1}-1.$ \\
Proceeding in this way, we get
$$deg(\alpha_{i_1}+ \alpha_{i_2}+\dots+\alpha_{i_l}+\beta_{j_1})=(q^l-1)q^{r-l}(q-1)q^{s-1}-1.$$
Now, run $\beta_{j_p}$ from $1\leq p\leq m$, by the above arguments, it is easy to see that
$$deg(\alpha_{i_1}+\alpha_{i_2}+\dots+ \alpha_{i_l}+\beta_{j_1}+\beta_{j_2}+\dots+\beta_{j_m})=(q^l-1)q^{r-l}(q^m-1)q^{s-m}-1.$$
Hence,
$$deg(\alpha_{i_1}+\alpha_{i_2}+\dots+ \alpha_{i_l}+\beta_{j_1}+\beta_{j_2}+\dots+\beta_{j_m})=(q^l-1)(q^m-1)q^{n-(l+m)}-1.$$\qed

\section{$\Gamma_{U\oplus W}(\mathbb{V})$ is not Eulerian but triangulated}

In this section, we first find the minimum degree and the edge connectivity of $\Gamma_{U\oplus W}(\mathbb{V})$. The main aim of this section is to show that the graph $\Gamma_{U\oplus W}(\mathbb{V})$ is not Eulerian but triangulated. To do this, we require the minimum degree of $\Gamma_{U\oplus W}(\mathbb{V})$, which is obtained as follows.
\begin{theorem}\label{T11} Let $\mathbb{V}=U\oplus W$ be an $n$ dimensional vector space over a field $\mathcal{F}$ with $q$ elements.  Then the minimum degree $\delta$ of $\Gamma_{U\oplus W}(\mathbb{V})$ is $(q-1)^2q^{n-2} -1$.
\end{theorem}
\noindent\textbf{Proof.} From Theorem \ref{T7}, the degree of the vertex $c_1\alpha_{i_1}+c_{2}\alpha_{i_2}+\dots+ c_{l}\alpha_{i_l}+d_{1}\beta_{j_1}+d_{2}\beta_{j_2}+\dots+ d_{m}\beta_{j_m}$, where $c_1c_2\dots c_l d_1d_2\dots d_m \neq 0$, is $(q^l-1)(q^m-1)q^{n-(l+m)}-1$. Thus the degree will be minimized if
$l =m= 1$  and hence $\delta = (q-1)^2q^{n-2} -1$. \qed
\begin{corollary}\label{C1} Let $\mathbb{V}=U\oplus W$ be an $n$ dimensional vector space over a field $\mathcal{F}$ with $q$ elements.  Then the edge connectivity  of $\Gamma_{U\oplus W}(\mathbb{V})$ is $(q-1)^2q^{n-2} -1$.
\end{corollary}
\noindent\textbf{Proof.}
 From Theorem \ref{T2}, if $dim(\mathbb{V})=2$, then $\Gamma_{U\oplus W}(\mathbb{V})$ is complete. Thus, the degree of every vertex is same. Therefore, the edge connectivity of $\Gamma_{U\oplus W}(\mathbb{V})$  is $(q-1)^2q^{n-2} -1$. On the other hand, if $dim(\mathbb{V})\geq3$, then  $diam(\Gamma_{U\oplus W}(\mathbb{V}))=2$, by Theorem \ref{T1}. Therefore, its edge connectivity is equal to its minimum degree (see \cite{Ple},). i.e, $(q-1)^2q^{n-2} -1$. \qed

Now, we show that the graph $\Gamma_{U\oplus W}(\mathbb{V})$ is not Eulerian.

\begin{theorem}\label{T13} Let $\mathbb{V}=U\oplus W$ be an $n$ dimensional vector space over a field $\mathcal{F}$ with $q$ elements.  Then   $\Gamma_{U\oplus W}(\mathbb{V})$ is not an Eulerian.
\end{theorem}
\noindent\textbf{Proof.}
\indent If $q=2$, by Theorem \ref{T7}, every vertex is of odd degree. Again, if $q$ is odd prime,  $\Gamma_{U\oplus W}(\mathbb{V})$ is not Eulerian, since, in this case the degree of every vertex is odd. Thus the graph is not Eulerian in any case.\qed

In order to determine whether the graph is triangulated or not, we need to find the order and size  of $\Gamma_{U\oplus W}(\mathbb{V})$ and  the  minimum degree of a vertex.
\begin{theorem}\label{T8}Let $\mathbb{V}=U\oplus W$ be an $n$ dimensional vector space over a field $\mathcal{F}$ with $q$ elements. Then the order of $\Gamma_{U\oplus W}(\mathbb{V})$ is $(q^r-1)(q^s-1)$ and size $M$ of $\Gamma_{W}(\mathbb{V}_\alpha)$ is$$\frac{\big(q^{2r}-(2q-1)^r\big)\big(q^{2s}-(2q-1)^s\big)-(q^r-1)(q^s-1)}{2}.$$
\end{theorem}
\noindent\textbf{Proof.} It is easy to observe that the order of $\Gamma_{U\oplus W}(\mathbb{V})$ is $(q^r-1)(q^s-1)$.\\
It is easy to see that the number of vectors having exactly $l$ $\alpha_i$'s and $m$, $\beta_j$'s with nonzero coefficient in its basic representation are  $\binom{r}{l}(q-1)^l\binom{s}{m}(q-1)^m$, where $1\leq l\leq r$ and $1\leq m\leq s$ .  By Theorem \ref{T7} and noting the fact that the sum of the degrees of all the vertices in  $\Gamma_{U\oplus W}(\mathbb{V})$ is equal to $2M$, we have
\begin{align*}
2M&= \sum_{l=1}^{r}\sum_{m=1}^{s}\binom{r}{l}(q-1)^l\binom{s}{m}(q-1)^m\bigg( (q^l-1)(q^m-1)q^{n-(l+m)}-1\bigg)\\
 &=\sum_{l=1}^{r}\binom{r}{l}(q-1)^l(q^l-1)q^{r-l}\sum_{m=1}^{s}\binom{s}{m}(q-1)^m (q^m-1)q^{s-m}\\
&-\sum_{l=1}^{r}\binom{r}{l}(q-1)^l\sum_{m=1}^{s}\binom{s}{m}(q-1)^m \\
&=\sum_{l=1}^{r}\binom{r}{l}(q-1)^l(q^r-q^{r-l})\sum_{m=1}^{s}\binom{s}{m}(q-1)^m (q^s-q^{s-m}\\
&-\sum_{l=1}^{r}\binom{r}{l}(q-1)^l\sum_{m=1}^{s}\binom{s}{m}(q-1)^m\\
\end{align*}
That is,
\begin{align*}
2M&=\Bigg(\sum_{l=1}^{r}\binom{r}{l}(q-1)^lq^r-\sum_{l=1}^{r}\binom{r}{l}(q-1)^lq^{r-l}\Bigg)\Bigg(\sum_{m=1}^{s}\binom{s}{m}(q-1)^m q^s\\&-\sum_{m=1}^{s}\binom{s}{m}(q-1)^mq^{s-m}\Bigg)
-\sum_{l=1}^{r}\binom{r}{l}(q-1)^l\sum_{m=1}^{s}\binom{s}{m}(q-1)^m\\
&= \big(q^r(q^r-1)-[(q+q-1)^r-q^r]\big)\big(q^s(q^s-1)-[(q+q-1)^s-q^s]\big)\\
&-[(q-1+1)^r-1][(q-1+1)^r-1]\\
&=\big(q^{2r}-(2q-1)^r\big)\big(q^{2s}-(2q-1)^s\big)-(q^r-1)(q^s-1).\\
\end{align*}
Hence the result follows.\qed

 The following theorem shows that $\Gamma_{U\oplus W}(\mathbb{V})$ is triangulated.

\begin{theorem}\label{T9} Let $\mathbb{V}=U\oplus W$ be an $n$ dimensional vector space over a field $\mathcal{F}$ with $q$ elements.
\item{(i)} If $n=2$, $dim(W)= n-1$ and $q=2$, then $\Gamma_{U\oplus W}(\mathbb{V})$ is a trivial graph.
\item{(ii)} If $n=2$, $dim(W)= n-1$ and $q\neq 2$, then $\Gamma_{U\oplus W}(\mathbb{V})$ is triangulated.
\item{(iii)} If $n=3$, $dim(W)= n-1$ and $q= 2$,  then $\Gamma_{U\oplus W}(\mathbb{V})$ contains no cycle.
\item{(iv)} If $n=3$, $dim(W)= n-1$ and $q\neq 2$, then $\Gamma_{U\oplus W}(\mathbb{V})$ is triangulated.
\item{(v)} If $n\geq 4$, $dim(W)= n-1$ and $q\geq 2$, then $\Gamma_{U\oplus W}(\mathbb{V})$ is triangulated.
\item{(vi)} If $n\geq 4$,  $dim(W)\leq n-2$ and $q\geq 2$, then $\Gamma_{U\oplus W}(\mathbb{V})$ is triangulated.
\end{theorem}
\noindent\textbf{Proof.}
\noindent (i). For $n=2$ and $q=2$, the vertex set $\mathbf{V}$ is $\{\alpha_1+\beta_1\}$. Thus,  $\Gamma_{U\oplus W}(\mathbb{V})$ is a  trivial graph.
\noindent (ii). For $n=2$ and $q\neq 2$,  $\Gamma_{U\oplus W}(\mathbb{V})$  is a complete graph and therefore triangulated.
\noindent (iii). For $n=3$ and $q= 2$, the vertex set $\mathbf{V}$ is $\{\alpha_1+\beta_1,\alpha_1+\beta_2,\alpha_1+\beta_1+\beta_2\}$. Clearly, $\Gamma_{U\oplus W}(\mathbb{V})$ contains no cycle.
\noindent (iv). For $n=3$ and $q\neq 2$,  there exists $a\in \mathcal{F}\setminus \{0, 1\}$. For any arbitrary $\alpha_1+\beta_1+\beta_2$, there exist two vertices either$\{ \alpha_1+a\beta_1+a\beta_2,\alpha_1+\beta_1\}$ or $\{ \alpha_1+a\beta_1+a\beta_2,\alpha_1+\beta_2\}$  such that $\alpha_1+\beta_1+\beta_2\sim \alpha_1+a\beta_1+a\beta_2\sim \alpha_1+\beta_1\sim \alpha_1+\beta_1+\beta_2$ and $\alpha_1+\beta_1+\beta_2\sim \alpha_1+a\beta_1+a\beta_2\sim \alpha_1+\beta_2\sim \alpha_1+\beta_1+\beta_2$ is a cycle of length three. Hence $\Gamma_{U\oplus W}(\mathbb{V})$  is  triangulated.
\noindent (v). For $n\geq 4$ and $q\geq 2$, every vertex is of the form $c_1\alpha_1+d_1\beta_1+d_2\beta_2+\dots+\dots+d_k\beta_k$, where $c_1\neq0,d_1,d_2,d_k\in \mathcal{F}$ and $1\leq k\leq n-1$. Consider an arbitrary vertex $x=c_1\alpha_1+d_1\beta_{i_1}$, there exists  $y=c_1\alpha_1+d_1\beta_{i_1}+d_2\beta_{i_2}$ and $z=c_1\alpha_1+d_1\beta_{i_1}+d_2\beta_{i_1}+d_3\beta_{i_1}$ such that $x\sim y\sim z\sim x$ is a cycle of length three. Hence $\Gamma_{U\oplus W}(\mathbb{V})$  is  triangulated.
\noindent (vi) For $n\geq4$ and $q\geq2$, every vertex is of the form  $c_1\alpha_{i_1}+c_2\alpha_{i_2}+\dots+c_l\alpha_{i_l}+d_1\beta_{j_1}+d_2\beta_{j_1}+\dots+d_m\beta_{j_m}$, where $l+m=n.$
Consider any arbitrary vertex $x=c_1\alpha_{i_1}+c_2\alpha_{i_2}+\dots+c_l\alpha_{i_l}+d_1\beta_{j_1}+d_2\beta_{j_1}+\dots+d_m\beta_{j_m}$. There exists  $y=c_p\alpha_{i_p}+c_q\alpha_{i_q}+d_p\beta_{j_p}+d_q\beta_{j_q}$ and $z=c_p\alpha_{i_p}+d_p\beta_{j_p}$ such that $x,y$ and $z$ share at least two basis vectors having nonzero coefficients in there basic representation. In other words $x\sim y\sim z\sim x$ is a cycle of length three. Hence $\Gamma_{U\oplus W}(\mathbb{V})$  is  triangulated.\qed
Now, we obtain the girth of $\Gamma_{U\oplus W}(\mathbb{V})$.
\begin{theorem}\label{T10} Let $\mathbb{V}=U\oplus W$ be an $n$ dimensional vector space over a field $\mathcal{F}$ with $q$ elements. Then girth
\[
gr(\Gamma_{W}(\mathbb{V}_\alpha))=
\begin{cases}
\infty& \quad \text{ if \it{n=2, dim(W)= n-1 and q=2 }}\\
3& \quad \text{if \it{n=2, dim(W)= n-1 and q$\neq$ 2}}\\
\infty& \quad \text{if \it{n=3, dim(W)= n-1 and q= 2}}\\
3& \quad \text{if \it{n=3, dim(W)= n-1 and q$\neq$ 2}}\\
3& \quad \text{if \it{n$\geq$ 4, dim(W)= n-1 and q$\geq$ 2}}\\
3& \quad \text{if \it{n$\geq$ 4,  dim(W)$\leq$ n-2 and q$\geq$ 2}}.\\
\end{cases}
\]
\end{theorem}
\noindent\textbf{Proof.}  The result follows from Theorem \ref{T9}. \qed

\section{Maximal Cliques in $\Gamma_{U\oplus W}(\mathbb{V})$}

In this section, we find the clique
number of $\Gamma_{U\oplus W}(\mathbb{V})$.
For $x\in \Gamma_{U\oplus W}(\mathbb{V})$, let $S_x$ (skeleton of $x$) be the set of $\alpha_i$’s and $\beta_j$'s with nonzero coefficients
in the basic representation of $x$ with respect to the basis of $\mathbb{V}$. It is to be noted that two distinct vertices may
have same skeleton. Moreover, if $x=c_1\alpha_{i_1}+c_2\alpha_{i_2}+\dots+c_{l-1}\alpha_{i_{l-1}}+d_1\beta_{j_1}+d_2\beta_{j_1}+\dots+d_{m-1}\beta_{j_{m-1}}$ and $y=c_1\alpha_{i_1}+c_2\alpha_{i_2}+\dots+c_l\alpha_{i_l}+d_1\beta_{j_1}+d_2\beta_{j_1}+\dots+d_m\beta_{j_m}$, then $S_x\subset S_y$.  Also, $2\leq |S_x |\leq n$, for all $x\in \Gamma_{U\oplus W}(\mathbb{V})$. Let $\mathbf{M}$ be a maximal clique in $ \Gamma_{U\oplus W}(\mathbb{V})$ and $S_x^U$ and $S_x^W$ be the skeleton of elements in subspaces $U$ and $W$ of $\mathbb{V}$. Define $S(\mathbf{M})=\{S_x=S_x^U\cup S_x^W: x\in \mathbf{M}\} $ and $S[\mathbf{M}]=\{|S_x|==S_x^U+S_x^W: S_x\in S(\mathbb{M})\}$. Since $S[\mathbf{M}]\neq \phi$, by well-ordering principle, it has a least element, say $k_1+k_2=k$.
Then, there exist some $x^\ast\in \mathbf{M}$ with $|S_x^\ast|=k_1+k_2$, where $x^\ast=c_1\alpha_{i_1}+c_2\alpha_{i_2}+\dots+c_{k_1}\alpha_{i_{k_1}}+d_1\beta_{j_1}+d_2\beta_{j_1}+\dots+d_{k_2}\beta_{j_{k_2}}$. Depending upon the choice of $k_1\leq \frac{r}{2}$, $k_2\leq \frac{s}{2}$ or $k_1> \frac{r}{2}$, $k_2> \frac{r}{2}$, we show that there exists four types of maximal cliques in $\Gamma_{U\oplus W}(\mathbb{V})$.

\begin{theorem}\label{T14} Let $\mathbf{M}$ be a maximal clique in $\Gamma_{U\oplus W}(\mathbb{V})$. If $k_1+k_2$ is the least element of $S[\mathbf{M}]$ with $k_1\leq \frac{r}{2}$ and $k_2\leq \frac{s}{2}$, then $\mathbf{M}\in M_{k_1,k_2}^{i,j}$ with $1\leq k_1\leq\frac{r}{2},~1\leq k_2\leq \frac{s}{2}$, $i\in\{1,2,\dots,r\}$  and $j\in\{1,2,\dots,r\}$, where $ M_{k_1,k_2}^{i,j}=\{x\in \Gamma_{U\oplus W}(\mathbb{V}):\alpha_{i},\beta_{j}\in S_x\text{and }     |S_x|\geq k\}$ and
\begin{equation*}
     |\mathbf{M}|=(q-1)^2\sum_{i=k_1-1}^{r-1}\sum_{j=k_2-1}^{s-1}\binom{r-1}{i}\binom{s-1}{j}(q-1)^{i+j}
\end{equation*}
\end{theorem}
\noindent\textbf{Proof.} Let $x=u+w\in \mathbf{M}$, where $u\in U$ and $w\in W$.  Since $k_1\leq \frac{r}{2}$, by Erdos-Ko-Rado theorem\cite{EKR}, the maximum number of pairwise-intersecting $k_1$-subsets in $U$ is $\binom{r-1}{k_1-1}$ and the maximum is achieved only if each $k_1$-subset contains a fixed element,
say $\alpha_{i}$.
Now, the number of $u$'s in $U$ with $|S_u|=k_1,k_1+1,k_1+2,\dots, r$ and $\alpha_{i}\in S_u$ are  \begin{eqnarray}\label{L16}
\binom{r-1}{k_1-1}(q-1)^{k_1}, \binom{r-1}{k_1}(q-1)^{k_1+1},\binom{r-1}{k_1+1}(q-1)^{k_1+2},\dots,\binom{r-1}{r-1}(q-1)^{r},
\end{eqnarray}    respectively.\\
Using the similar argument, for $k_2\leq \frac{s}{2}$, the number of $W$'s in $W$ with $|S_w|=k_2,k_2+1,k_2+2,\dots, s$ and $\beta_{j}\in S_w$ are  \begin{eqnarray}\label{L17}
\binom{s-1}{k_2-1}(q-1)^{k_2}, \binom{s-1}{k_2}(q-1)^{k_2+1},\binom{s-1}{k_2+1}(q-1)^{k_2+2},\dots,\binom{s-1}{s-1}(q-1)^{s},
\end{eqnarray}
respectively. As $\mathbf{M}$ is a maximal clique, and minimum of $S_x$, for $x\in \mathbf{M}$ is $k_1+k_2$ and $S_x\cap S_y\neq \phi$, so
$ \mathbf{M}=\{x\in \Gamma_{U\oplus W}(\mathbb{V}):\alpha_{i},\beta_{j}\in S_x\text{and }     |S_x|\geq k_1+k_2\}$.\\
Now, the number of $x$'s in $\mathbf{M}$ with $|S_x|=k_1+k_2$ and $\alpha_{i},\beta_{j}\in S_x$ is $$\binom{r-1}{k_1-1}(q-1)^{k_1}\binom{s-1}{k_2-1}(q-1)^{k_2}.$$
The number of $x$'s in $\mathbf{M}$ with $|S_x|=k_1+k_2+1$ and $\alpha_{i},\beta_{j}\in S_x$ are $$ \binom{r-1}{k_1-1}(q-1)^{k_1}\binom{s-1}{k_2}(q-1)^{k_2+1}+ \binom{r-1}{k_1}(q-1)^{k_1+1} +\binom{s-1}{k_2-1}(q-1)^{k_2}.$$
Similarly, for $|S_x|= k_1+k_2+2,\dots, r+s$ and $\alpha_{i},\beta_{j}\in S_x$, we can find the number of $x$'s in $\mathbf{M}$ by the help of equations (\ref{L16}) and (\ref{L17}).
 Therefore, we have
\begin{eqnarray*}
|\mathbf{M}|&=&\binom{r-1}{k_1-1}(q-1)^{k_1}\binom{s-1}{k_2-1}(q-1)^{k_2}+
\binom{r-1}{k_1-1}(q-1)^{k_1}\binom{s-1}{k_2}(q-1)^{k_2+1}\\
&+& \binom{r-1}{k_1}(q-1)^{k_1+1} +\binom{s-1}{k_2-1}(q-1)^{k_2}+\dots+ \binom{r-1}{r-1}(q-1)^{r}\binom{s-1}{s-1}(q-1)^{s}\\
&=& \binom{r-1}{k_1-1}(q-1)^{k_1} \Bigg[\binom{s-1}{k_2-1}(q-1)^{k_2}+ \binom{s-1}{k_2}(q-1)^{k_2+1}\\
&&+\dots+ \binom{s-1}{s-1}(q-1)^{s}\Bigg]+\binom{r-1}{k_1}(q-1)^{k_1+1}\Bigg[\binom{s-1}{k_2-1}(q-1)^{k_2}\\
&&+ \binom{s-1}{k_2}(q-1)^{k_2+1}
+\dots+ \binom{s-1}{s-1}(q-1)^{s}\Bigg]+\dots\\
&&+\binom{r-1}{r-1}(q-1)^{r}\Bigg[\binom{s-1}{k_2-1}(q-1)^{k_2}+ \binom{s-1}{k_2}(q-1)^{k_2+1}\\
&&+\dots+ \binom{s-1}{s-1}(q-1)^{s}\Bigg]\\
&=& \sum_{i=k_1-1}^{r-1}\binom{r-1}{i}(q-1)^{i+1}\sum_{j=k_2-1}^{s-1}\binom{s-1}{j}(q-1)^{j+1}\\
&=& (q-1)^2\sum_{i=k_1-1}^{r-1}\sum_{j=k_2-1}^{s-1}\binom{r-1}{i}\binom{s-1}{j}(q-1)^{i+j}.
\end{eqnarray*}
It is to be noted that for same value of $k_1,k_2$ and by fixing different $\alpha_i$'s and $\beta_j$'s, we get different
maximal cliques. Since these maximal cliques depends both on $k_1,k_2$ and  $\alpha_i, \beta_j$, we get a family
of maximal cliques $M_{k_1,k_2}^{i,j}=\{x\in \Gamma_{U\oplus W}(\mathbb{V}):\alpha_{i},\beta_{j}\in S_x\text{and } |S_x|\geq k\}$ , where $1\leq k_1\leq\frac{r}{2},1\leq k_2\leq \frac{s}{2}$, $i\in\{1,2,\dots,r\}$  and $j\in\{1,2,\dots,r\}$ and $\mathbf{M}\in M_{k_1,k_2}^{i,j}.$\qed

\begin{theorem}\label{T15} Let $\mathbf{M}$ be a maximal clique in $\Gamma_{U\oplus W}(\mathbb{V})$. If $k_1+k_2$ is the least element of $S[\mathbf{M}]$ with $k_1\leq \frac{r}{2}$ and $k_2> \frac{s}{2}$, then $\mathbf{M}\in M_{k_1,k_2}^{i,j}$ with $1\leq k_1\leq\frac{r}{2}$,$k_2=\lfloor\frac{s}{2}\rfloor+1$, $i\in\{1,2,\dots,r\}$  and $j\in\{1,2,\dots,r\}$,  where $ M_{k_1,k_2}^{i,j}=\{x\in \Gamma_{U\oplus W}(\mathbb{V}):\alpha_{i},\beta_{j}\in S_x \text{and}     |S_x|=|S_u|+|S_w|\geq k\}$ and
\begin{equation*}
     |\mathbf{M}|=(q-1)\sum_{i=k_1-1}^{r-1}\sum_{j=k_2}^{s}\binom{r-1}{i}\binom{s}{j}(q-1)^{i+j}.
\end{equation*}
\end{theorem}
\noindent\textbf{Proof.} Let $x=u+w\in \mathbf{M}$, where $u\in U$ and $w\in W$.  Since $k_1\leq \frac{r}{2}$, following the same arguments as in Theorem \ref{14}, the number of $u$'s in $U$ with $|S_u|=k_1,k_1+1,k_1+2,\dots, r$ and $\alpha_{i}\in S_u$ are  \begin{eqnarray*}
\binom{r-1}{k_1-1}(q-1)^{k_1}, \binom{r-1}{k_1}(q-1)^{k_1+1},\binom{r-1}{k_1+1}(q-1)^{k_1+2},\dots,\binom{r-1}{r-1}(q-1)^{r},
\end{eqnarray*}    respectively. Now, for $k_2> \frac{s}{2}$,  by Erdos-Ko-Rado theorem\cite{EKR}, the maximum number of pairwise-intersecting $k_2$-subsets in $V$ is $\binom{s}{k_2}$ and the maximum is achieved only if each $k_2$-subset contains a fixed element, say $\beta_{j}$.  The number of $w$'s in $W$ with $|S_w|=k_2,k_2+1,k_2+2,\dots, s$ and $\beta_{j}\in S_w$ are  \begin{eqnarray*}
\binom{s}{k_2}(q-1)^{k_2}, \binom{s}{k_2+1}(q-1)^{k_2+1},\binom{s}{k_2+2}(q-1)^{k_2+2},\dots,\binom{s}{s}(q-1)^{s},
\end{eqnarray*}
respectively. As $\mathbf{M}$ is a maximal clique,  and minimum of $S_x$ for $x\in \mathbf{M}$ is $k_1+k_2$ and $S_x\cap S_y\neq \phi$, so
$ \mathbf{M}=\{x\in \Gamma_{U\oplus W}(\mathbb{V}):\alpha_{i},\beta_{j}\in S_x\text{and }     |S_x|\geq k_1+k_2\}$.\\
Now, following the same procedure as in the proof of Theorem \ref{T14}, we have
\begin{align*}
|\mathbf{M}|&=\binom{r-1}{k_1-1}(q-1)^{k_1}\binom{s}{k_2}(q-1)^{k_2}+
\binom{r-1}{k_1-1}(q-1)^{k_1}\binom{s}{k_2+1}(q-1)^{k_2+1}\\
&+ \binom{r-1}{k_1}(q-1)^{k_1+1} +\binom{s}{k_2}(q-1)^{k_2}+\dots+ \binom{r-1}{r-1}(q-1)^{r}\binom{s}{s}(q-1)^{s}\\
&= \binom{r-1}{k_1-1}(q-1)^{k_1} \Bigg[\binom{s}{k_2}(q-1)^{k_2}+ \binom{s}{k_2+1}(q-1)^{k_2+1}\\
&+\dots+ \binom{s}{s}(q-1)^{s}\Bigg]+\binom{r-1}{k_1}(q-1)^{k_1+1}\Bigg[\binom{s}{k_2}(q-1)^{k_2}\\
&+ \binom{s}{k_2+1}(q-1)^{k_2+1}
+\dots+ \binom{s}{s}(q-1)^{s}\Bigg]+\dots
\end{align*}
\begin{align*}
&+\binom{r-1}{r-1}(q-1)^{r}\Bigg[\binom{s}{k_2}(q-1)^{k_2}+ \binom{s}{k_2+1}(q-1)^{k_2+1}\\
&+\dots+ \binom{s}{s}(q-1)^{s}\Bigg]\\
&= \sum_{i=k_1-1}^{r-1}\binom{r-1}{i}(q-1)^{i+1}\sum_{j=k_2}^{s}\binom{s}{j}(q-1)^{j}\\
&= (q-1)\sum_{i=k_1-1}^{r-1}\sum_{j=k_2}^{s}\binom{r-1}{i}\binom{s}{j}(q-1)^{i+j}.
\end{align*}
It is to be noted that for same value of $k_1,k_2$ and by fixing different $\alpha_i$'s and $\beta_j$'s, we get different
maximal cliques. Since these maximal cliques depend both on $k_1,k_2$ and  $\alpha_i, \beta_j$, we get a family
of maximal cliques $M_{k_1,k_2}^{i,j}=\{x\in \Gamma_{U\oplus W}(\mathbb{V}):\alpha_{i},\beta_{j}\in S_x\text{and } |S_x|=|S_u|+|S_w|\geq k_1+k_2\},$ where $1\leq k_1\leq\frac{r}{2},k_2> \frac{s}{2}$, $i\in\{1,2,\dots,r\}$  and $j\in\{1,2,\dots,r\}$ and $\mathbf{M}\in M_{k_1,k_2}^{i,j}.$ Now, as $k_2>\frac{s}{2}$,  $\{w\in W:  |S_w|\geq k_2+1\}\subset \{w\in W:  |S_w|\geq k_2\}$. Maximality of $\mathbf{M}$ is obtained by minimizing $k_2$ provided $k_2>\frac{s}{2}$, that is $k_2=\lfloor\frac{s}{2}\rfloor+1$. Thus, we get a family of maximal cliques $M_{k_1}^{i,j}$,  where $1\leq k_1\leq\frac{r}{2}$, $i\in\{1,2,\dots,r\}$  and $j\in\{1,2,\dots,r\}$ and $\mathbf{M}\in M_{k_1}^{i,j}.$
\qed

\begin{theorem}\label{T16} Let $\mathbf{M}$ be a maximal clique in $\Gamma_{U\oplus W}(\mathbb{V})$. If $k_1+k_2$ is the least element of $S[\mathbf{M}]$ with $k_1> \frac{r}{2}$ and $k_2\leq \frac{s}{2}$, then $\mathbf{M}\in M_{k_1,k_2}^{i,j}$ with $k_1=\lfloor\frac{r}{2}\rfloor+1$, $1\leq k_2\leq\frac{s}{2}$, $i\in\{1,2,\dots,r\}$  and $j\in\{1,2,\dots,r\}$,  where $ M_{k_1,k_2}^{i,j}=\{x\in \Gamma_{U\oplus W}(\mathbb{V}):\alpha_{i},\beta_{j}\in S_x \text{and}     |S_x|=|S_u|+|S_w|\geq k\}$ and
\begin{equation*}
     |\mathbf{M}|=(q-1)\sum_{i=k_1}^{r}\sum_{j=k_2-1}^{s-1}\binom{r}{i}\binom{s-1}{j}(q-1)^{i+j}.
\end{equation*}
\end{theorem}
\noindent\textbf{Proof.} The proof is similar to that of Theorems \ref{T14} and \ref{T15}.\qed

\begin{theorem}\label{T17} Let $\mathbf{M}$ be a maximal clique in $\Gamma_{U\oplus W}(\mathbb{V})$. If $k_1+k_2$ is the least element of $S[\mathbf{M}]$ with $k_1> \frac{r}{2}$ and $k_2> \frac{s}{2}$, then $k_1=\lfloor\frac{r}{2}\rfloor+1$ and  $k_2=\lfloor\frac{s}{2}\rfloor+1$ and $\mathbf{M}=\{x\in \Gamma_{U\oplus W}(\mathbb{V}): |S_x|=|S_u|+|S_w|,~ |S_u| \geq \lfloor\frac{r}{2}\rfloor+1 ~\text{and}~ |S_w| \geq \lfloor\frac{s}{2}\rfloor+1\}$ and
\begin{equation*}
     |\mathbf{M}|=\sum_{i=\lfloor\frac{r}{2}\rfloor+1}^{r}\sum_{j=\lfloor\frac{s}{2}\rfloor+1}^{s}\binom{r}{i}\binom{s}{j}(q-1)^{i+j}.
\end{equation*}
\end{theorem}
\noindent\textbf{Proof.} Let $x=u+w\in \mathbf{M}$, where $u\in U$ and $w\in W$.  Since $k_1> \frac{r}{2}$ and $k_2> \frac{s}{2}$, by following the same arguments as in Theorems \ref{T14} and \ref{T15}, the number of $u$'s in $U$ with $|S_u|=k_1,k_1+1,k_1+2,\dots, r$ and $\alpha_{i}\in S_u$ are
\begin{eqnarray*}
\binom{r}{k_1}(q-1)^{k_1}, \binom{r}{k_1+1}(q-1)^{k_1+1},\binom{r}{k_1+2}(q-1)^{k_1+2},\dots,\binom{r}{r}(q-1)^{r},
\end{eqnarray*}
Similarly, the number of $W$'s in $W$ with $|S_w|=k_2,k_2+1,k_2+2,\dots, s$ and $\beta_{j}\in S_w$ are
\begin{eqnarray*}
\binom{s}{k_2}(q-1)^{k_2}, \binom{s}{k_2+1}(q-1)^{k_2+1},\binom{s}{k_2+2}(q-1)^{k_2+2},\dots,\binom{s}{s}(q-1)^{s},
\end{eqnarray*}
respectively. As $\mathbf{M}$ is a maximal clique,  and minimum of $S_x$ for $x\in \mathbf{M}$ is $k_1+k_2$ and $S_x\cap S_y\neq \phi$, so
$ \mathbf{M}=\{x\in \Gamma_{U\oplus W}(\mathbb{V}):\alpha_{i},\beta_{j}\in S_x\text{and }     |S_x|\geq k_1+k_2\}$.\\
Therefore, following the same procedure as in the proof of Theorem \ref{T14}, we have
\begin{eqnarray*}
|\mathbf{M}|&=&\sum_{i=k_1}^{r}\binom{r}{i}(q-1)^{i}\sum_{j=k_2}^{s}\binom{s}{j}(q-1)^{j}\\
&=& \sum_{i=k_1}^{r}\sum_{j=k_2}^{s}\binom{r}{i}\binom{s}{j}(q-1)^{i+j}.
\end{eqnarray*}
 Now, as $k_1>\frac{r}{2}$ and $k_2>\frac{s}{2}$,  $\{x\in \Gamma_{U\oplus W}(\mathbb{V}):  |S_x|\geq k_1+k_2+1\}\subset \{x\in \Gamma_{U\oplus W}(\mathbb{V}):  |S_x|\geq k_1+k_2\}$. Maximality of $\mathbf{M}$ is obtained by minimizing $k_1$ and $k_2$ provided $k_1>\frac{r}{2}$ and $k_2>\frac{s}{2}$, that is, $k_1=\lfloor\frac{r}{2}\rfloor+1$ and  $k_2=\lfloor\frac{s}{2}\rfloor+1$. Thus,
 \begin{eqnarray*}
|\mathbf{M}|&=& \sum_{i=\lfloor\frac{r}{2}\rfloor+1}^{r}\sum_{j=\lfloor\frac{s}{2}\rfloor+1}^{s}\binom{r}{i}\binom{s}{j}(q-1)^{i+j}.
\end{eqnarray*}
\qed

\begin{remark}\label{R1}
Note that $|M_{k_1,k_2}^{i,j}|$ is maximum when $k_1=k_2=1$ in case $k_1\leq \frac{r}{2}$ and $k_2\leq \frac{s}{2}$. So, combining Theorems \ref{T14}, \ref{T15}, \ref{T16} and \ref{T17}, it is easy to see that the clique number of $\Gamma_{U\oplus W}(\mathbb{V})$ is
\begin{eqnarray*}
\omega(\Gamma_{W}(\mathbb{V}_\alpha))=max\Bigg\{(q-1)q^{r-1}(q-1)q^{s-1}, (q-1)q^{r-1}\sum_{j=\lfloor\frac{s}{2}\rfloor+1}^{s}\binom{s}{j}(q-1)^{j},\\ \sum_{i=\lfloor\frac{r}{2}\rfloor+1}^{r}\binom{r}{i}(q-1)^{i}(q-1)q^{s-1}, \sum_{i=\lfloor\frac{r}{2}\rfloor+1}^{r}\binom{r}{i}(q-1)^{i}\sum_{j=\lfloor\frac{s}{2}\rfloor+1}^{s}\binom{s}{j}(q-1)^{j}.
\Bigg\}
\end{eqnarray*}
 and it depends on $q$, $r$, $s$.
 \end{remark}

\begin{corollary}\label{C2} If $\mathcal{F}= \mathbb{F}_2$ and $n\geq 4$, then the clique number $\omega(\Gamma_{U\oplus W}(\mathbb{V}))=2^{n-2}.$
\end{corollary}
\noindent\textbf{Proof.} If $q=2$, then $(q-1)q^{r-1}=q^{r-1}$. Moreover, it is easy to see that for $q=2$ and $p=2m$ or $p=2m+1$, $$\sum_{j=\lfloor\frac{p}{2}\rfloor+1}^{p}\binom{p}{j}(q-1)^{j}<2^{p-1}$$.

Using Remark \ref{R1}, we get $\omega(\Gamma_{U\oplus W}(\mathbb{V}))=2^{n-2}.$ \qed

\begin{corollary}\label{C3} If $\mathcal{F}= \mathbb{F}_2$, $n\geq 4$, and $\chi(\Gamma_{U\oplus W}(\mathbb{V}))$ is the chromatic number of $\Gamma_{U\oplus W}(\mathbb{V})$, then $$ 2^{n-2}\leq \chi(\Gamma_{U\oplus W}(\mathbb{V}))\leq \frac{q^n-(q^r+q^s+rs)}{2}+2^{n-3}+1 $$
\end{corollary}
\noindent\textbf{Proof.} For any graph $G$, $\omega(G)\leq \chi(G).$ Therefore, first part of inequality follows from corollary \ref{C2}. For the other inequality,  use the following result of  \cite{BrDu}
$$\chi(G)\leq \frac{\omega(G)+|\mathbf{V}|+1-\alpha(\Gamma_{W})}{2},$$ where $\alpha(\Gamma_{W})$ is the independence number of $G$. Thus, by using Theorems \ref{T2} and \ref{T10}, we have
\begin{eqnarray}\nonumber
\chi(\Gamma_{W}(\mathbb{V}_\alpha))&\leq&\frac{2^{n-2}+(q^r-1)(q^s-1)+1-rs}{2}\\
\nonumber &=&\frac{q^n-(q^r+q^s+rs)}{2}+2^{n-3}+1 .
\end{eqnarray}\qed

\section{Conclusion}

In this work, we have introduced the direct sum graph $\Gamma_{U\oplus W}(\mathbb{V})$ on a finite dimensional vector space $\mathbb{V}$. We investigated  some  basic  properties  like  connectedness,  completeness, independence number and domination number of $\Gamma_{U\oplus W}(\mathbb{V})$.  For $n\geq3$ and $q> 2$, we have shown that $\Gamma_{U\oplus W}(\mathbb{V})$ is triangulated. Also, we also provided the exact value of the minimum degree and the edge connectivity of $\Gamma_{U\oplus W}(\mathbb{V}))$. Moreover, we have seen that $\Gamma_{U\oplus W}(\mathbb{V})$ is not Eulerian. Further, we have found the clique number and the chromatic number of  $\Gamma_{U\oplus W}(\mathbb{V})$. For future research, we can think of investigating the following. (i) When $\Gamma_{U\oplus W}(\mathbb{V})$ is regular, (ii) whether $\Gamma_{U\oplus W}(\mathbb{V})$ is a line graph, (iii) characterize line graphs of $\Gamma_{U\oplus W}(\mathbb{V})$, (iv) to find genus of $\Gamma_{U\oplus W}(\mathbb{V})$.\\

\noindent\textbf{Acknowledgement} The research of S. Pirzada is supported by the National Board for Higher Mathematics (NBHM) research project number NBHM/02011/20/2022. Also the research of Dr. Bilal Ahmad Wani is supported by Dr. D.S. Kothari Post-Doctoral Fellowship Scheme Award Letter No. F.4-2/2006 (BSR)/MA/18-19/0037.\\

\noindent{\bf Data Availibility} Data sharing is not applicable to this article as no datasets were generated or analyzed
during the current study.\\

\noindent{\bf Conflict of Interest} The authors declared that they have no conflict of interest.

\end{document}